**On Validity of Reed's Conjecture for Classes of Graphs with Two Forbidden Subgraphs**

Medha Dhurandhar

**Abstract:** Reed's Conjecture is open for more than 20 years now. Here we prove that Reed's Conjecture is valid for (1) $\{P_4 \cup K_1$, Kite$\}$-free graphs (2) {Chair, Kite}-free graphs (3) $\{K_2 \cup \overline{K_2}$, H$\}$-free graphs and (4) $\{2K_2$, M$\}$-free graphs where H and M are graphs on six vertices each. Reed's conjecture is still open in general.

**Introduction:**

One of the most prominent problems in combinatorial optimization is to decide whether the vertices of a graph can be feasibly colored by a fixed number of different colors. If this fixed number is at least three, the aforementioned decision problem is known to be NP-complete. The associated optimization problem consists of computing the chromatic number $\chi$ of a graph. A lot of work is done to decide optimal bounds for the chromatic number of a graph. A lower bound for $\chi$ is the clique number $\omega$. A classical upper bound for $\chi$ in terms of the maximum degree $\Delta$ is provided by Brooks' Theorem, which states that $\chi \leq \Delta + 1$.

We consider here simple and undirected graphs. For terms which are not defined herein we refer to Bondy and Murty [1]. In 1998, Reed [2] proposed the following conjecture which gives, for any graph G, an upper bound for its chromatic number $\chi(G)$ in terms of the clique number $\omega(G)$ and the maximum degree $\Delta(G)$.

**Reed's Conjecture []:** *For any graph* G, $\chi(G) \leq \lceil \frac{\Delta + \omega + 1}{2} \rceil$.

In [3], [4], [5], [6], [7], [8], [9], [10], [11], [12], [13], [14], [15] it is shown that Reed's Conjecture holds for some graph classes defined by forbidden configurations:

- $\{P_5$, Flag$^C)$-free graphs,
- ($P_5$, $P_2 \cup P_3$, House, Dart)-free graphs,
- ($P_5$, Kite, Bull, $(K_3 \cup K_1) + K_1$)-free graphs,
- ($P_5$, C4)-free graphs,
- (Chair, House, Bull, $K_1 + C_4$)-free graphs,
- (Chair, House, Bull, Dart)-free graphs.
- $3K_1$-free graphs
- $\{2K_2, C_4\}$-free graphs
- Quasiline graphs
- $K_{1,3}$-free
- Generalized line graphs
- Graphs with $\chi \leq \omega + 2$
- Planar and toroidal graphs
- Decomposable graphs
- Perfect graphs
- Line graphs of Multigraphs
- Graphs with disconnected complements
- Graphs G with $\chi(G) > \lceil \frac{V(G)}{2} \rceil$ and graphs G with $\Delta(G) > \lceil \frac{V(G) - \alpha(G) + 3}{2} \rceil$
- Graphs G with $\Delta(G) \geq |V(G)| - 7$, and graphs G with $\Delta(G) \geq |V(G)| - \alpha(G) - 4$

This paper proves that Reed's Conjecture holds for

1. $\{P_4 \cup K_1$, Kite$\}$-free graphs,
2. {Chair, Kite}-free graphs,
3. $\{K_2 \cup \overline{K_2}$, H$\}$-free graphs and
4. $\{2K_2$, M$\}$-free graphs

where H and M are shown in **Figure 1**. Thus we prove validity of this conjecture for graphs with only two forbidden subgraphs.

H: 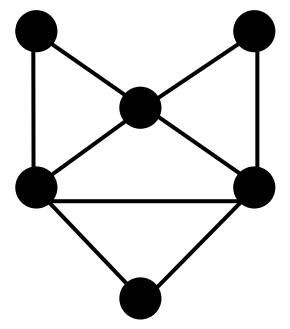 ; M: 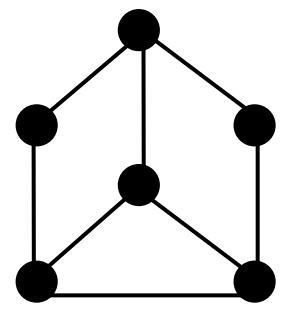

**Figure 1**

**Notation:** For a graph G, V(G), E(G), $\Delta$(G), $\omega$(G), $\chi$(G) denote the vertex set, edge set, maximum degree, size of a maximum clique, chromatic number respectively. For u $\in$ V(G), N(u) = {v $\in$ V(G) / uv $\in$ E(G)}, and $\overline{N(u)}$ = N(u)$\cup$(u). If S $\subseteq$ V(G), then <S> denotes the subgraph of G induced by S. If C is some coloring of G and if a vertex u of G is colored m in C, then u is called a m-vertex.

**Theorem 1:** If G is a {$P_4\cup K_1$, Kite}-free, then $\chi(G) \le \lceil \frac{\Delta+\omega+1}{2} \rceil$.

Proof: Let G be a smallest {$P_4\cup K_1$, Kite}-free graph with $\chi(G) > \lceil \frac{\Delta+\omega+1}{2} \rceil$. Let u $\in$ V(G). By minimality, $\chi(G) - 1 \le \chi(G-u) \le \lceil \frac{\Delta(G-u)+\omega(G-u)+1}{2} \rceil \le \lceil \frac{\Delta+\omega+1}{2} \rceil < \chi(G)$. Thus $\chi(G-u) = \chi(G) - 1 = \lceil \frac{\Delta+\omega+1}{2} \rceil$. Let C be a $\chi$(G)-coloring of G, in which only u is colored $\chi$ and N(u) has r vertices with unique colors in N(u). Let R = {X $\in$ N(u)/ X receives a unique color in N(u)}. Then |R| = r. Also $\Delta \ge \deg u \ge r + 2(\lceil \frac{\Delta+\omega+1}{2} \rceil - r)$ and r $\ge$ $\omega$(G)+1. Let Q $\subseteq$ R be s.t. <Q> is a maximum clique in R. Then |R-Q| $\ge$ 2 as r $\ge$ $\omega$(G)+1.

Note that as G is $P_4\cup K_1$-free, every chordless path of G is $P_4$ or $P_5$.

We have

**A. If $A_i$, $A_i$, $A_k \in$ R with $A_iA_j$, $A_iA_k \notin$ E(G) and {$A_i$, $a_j$, $a_i$, $A_j$} is an i-j bi-color path, then $A_ka_i \in$ E(G)** (else if $A_ka_i \notin$ E(G) and $A_ka_i' \in$ E(G), then <$A_i$, u, $A_k$, $a_i$', $a_i$> = $P_4\cup K_1$).

**B. Every vertex of R is non-adjacent to at most one vertex of R.**

Let if possible $A_iA_j$, $A_iA_k \notin$ E(G) where $A_i$, $A_j$, $A_k \in$ R. Let {$A_i$, $a_j$, $a_i$, $A_j$} be an i-j bi-color path**,** then by **A,** $A_ka_i \in$ E(G). Now $A_jA_k \notin$ E(G) (else <$a_i$, $A_k$, $A_j$, u, $A_i$> = Kite). Then by **A**, $A_ka_j \in$ E(G). Let $a_k$ be a k-vertex of $A_i$ on an i-k bi-color path from $A_i$ to $A_k$. Then by **A**, $A_ja_k \in$ E(G). Also either $a_ia_k \notin$ E(G) or $a_ja_k \notin$ E(G) (else <$a_i$, $a_j$, $a_k$, $A_i$, u> = Kite). W.l.g. let $a_ja_k \notin$ E(G). Then $A_j$ ($A_k$) is a unique j-vertex (k-vertex) of $a_k$ ($a_j$) (else if $a_j$' is another j-vertex $a_k$, then <u, $A_j$, $a_k$, $a_j$', $a_j$> = $P_4\cup K_1$). $\Rightarrow$ $\exists$ a bi-color path P = {$A_j$, $a_k$', $a_j$', $A_k$} where $a_l \neq a_l$' for l = j, k and <$a_j$', $a_k$', $A_j$, u, $a_j$> = $P_4\cup K_1$, a contradiction.

Now <R-Q> is complete (else if $\exists$ say $A_1$, $A_2 \in$ R-Q s.t. $A_1A_2 \notin$ E(G), then by **B**, $A_1A_j \in$ E(G) $\forall$ $A_j \in$ Q). Let $A_1$, $A_2 \in$ R-Q and $A_3$, $A_4 \in$ Q be s.t. $A_1A_3$, $A_2A_4 \notin$ E(G). Let $T_i$ = {$a_i$ / $a_i$ is an i-vertex of $A_{i+2}$}, $1 \le i \le 2$ and $T_j$ = {$a_{j-2}$ / $a_{j-2}$ is an (j-2)-vertex of $A_j$}, $3 \le j \le 4$. Then $a_ka_i \in$ E(G) for i = 1, 2 and k = i+2 (else <u, $A_k$, $a_i$, $a_k$', $a_k$> = $P_4\cup K_1$). Also $a_iA_j \in$ E(G) $\forall$ $1 \le i, j \le 4$, $j \neq i$ (else <$A_i$, u, $A_j$, $A_{i+2}$, $a_i$> = Kite). Now $a_1$ is not adjacent to both $a_2$ and $a_4$ (else <$a_1$, $a_2$, $a_4$, $A_4$, u> = Kite). Again $a_1$ is not non-adjacent to both $a_2$ and $a_4$ (else <$A_1$, $a_2$, $a_4$, $A_3$, $a_1$> = Kite). Hence w.l.g. let $a_1a_2 \in$ E(G) and $a_1a_4 \notin$ E(G). Similarly it can be seen that $a_3a_4 \in$ E(G) and $a_2a_3 \notin$ E(G). Also $A_4$ is the only 4-vertex of every vertex of $T_1$ (else if x is a 4-vertex of some $a_1$ in $T_1$, then $x \neq a_4$ and <x, $a_1$, $A_4$, u, $a_4$> = $P_4 \cup K_1$). Similarly $A_3$ is the only 3-vertex of every vertex of $T_2$. Then color vertices of $T_1$ by 4, $A_4$ by 2, vertices of $T_2$ by 3, $A_3$ by 1, and u by 4, a contradiction.

This proves **Theorem 1**.

**Theorem 2:** If G is a {Chair, Kite}-free, then $\chi(G) \leq \lceil \frac{\Delta+\omega+1}{2} \rceil$.

Proof: Let G be a smallest {Chair, Kite}-free graph with $\chi(G) > \lceil \frac{\Delta+\omega+1}{2} \rceil$. Let $u \in V(G)$. By minimality, $\chi(G) - 1 \leq \chi(G-u) \leq \lceil \ \rceil \leq \lceil \frac{\Delta+\omega+1}{2} \rceil < \chi(G)$. Thus $\chi(G-u) = \chi(G)-1 = \lceil \frac{\Delta+\omega+1}{2} \rceil$. Let C be a $\chi(G)$-coloring of G, in which only u is colored $\chi$. Let R = {X ∈ N(u)/ X receives a unique color in N(u)} and |R| = r. Then $\Delta \geq \deg u \geq r + 2(\lceil \frac{\Delta+\omega+1}{2} \rceil - r)$ and $r \geq \omega(G)+1$. Also let $Q \subseteq R$ s.t. <Q> is a maximum clique in <R>. As $r \geq \omega(G)+1$, $|R-Q| \geq 2$. **I**

**Case 1: ∃ A, B ∈ R-Q, s.t. AB ∉ E(G).**

**Case 1.1: ∃ C ∈ Q s.t.CA, CB ∉ E(G).**
Let A, B, C be colored 1, 2, 3 resply. Clearly ∃ a 3-1 bi-color path P from C to A. Let Cd, Ae ∈ P. As G is Chair-free, d is a unique 1-vertex of C (else if d' is another 1-vertex of C, then <C, d, d', u, A> = Chair). Similarly e is a unique 3-vertex of A. Also dB ∈ E(G) (else <u, A, B, C, d> = Chair). Similarly eB ∈ E(G). Now de ∈ E(G) (else if f (≠ C) is another 3-vertex of d on P, then if fB ∈ E(G), <B, f, e, u, C> = Chair & if fB ∉ E(G), <d, f, C, B, e> = Chair). Similarly if f is a 2-vertex of C, then fA, fd, fe ∈ E(G). But then <d, f, e, A, u> = Kite, a contradiction.

**Case 1.2: Case 1.1 does not hold.**
Then |Q| ≥ 2 and ∃ C, D ∈ Q s.t. AD, BC ∉ E(G) and AC, BD ∈ E(G). Let A, B, C, D be colored 1, 2, 3, 4 resply. Let e be a 2-vertex of A. Then eC ∉ E(G) (else <e, A, C, u, B> = Kite) and De ∈ E(G) (else <D, C, u, A, e> = Kite). Also e is the only 2-vertex of A (else if x is another 2-vertex of A, then <A, x, e, u, B> = Chair). Let f be a 2-vertex of C. Then if Df ∈ E(G), <D, e, f, B, u> = Chair and if Df ∉ E(G), <D, e, B, C, f> = Chair, a contradiction.

**Case 2: <R-Q> is complete.**
Let A, B ∈ R-Q and C, D ∈ Q be s.t. AC, BD ∉ E(G). Let A, B, C, D be colored 1, 2, 3, 4 resply.

**Case 2.1: AD (BC) ∈ E(G) and BC (AD) ∉ E(G).**
Since the configuration <u, A, B, C, D> is same as in **Case 1.2**, proof is similar.

**Case 2.2: AD, BC ∈ E(G).**
Let e be a 3-vertex of A, f a 4-vertex of B, g a 1-vertex of C and h a 2-vertex of D. Then as G is Chair-free e, f, g, h are unique 3, 4, 1, 2 vertices of A, B, C, D resply. Also eB ∈ E(G) (else <C, u, B, A, e> = Kite) and eD ∈ E(G) (else <C, D, u, A, e> = Kite). Similarly fC, fA, gB, gD, hC, hA ∈ E(G).

Now eg ∈ E(G) (else let j be the 1-vertex of e on a 1-3 path from A to C. Then if Dj ∈ E(G), <j, e, A, D, g> = Kite and if Dj ∉ E(G), <D, A, u, e, j> = Kite). Similarly fh ∈ E(G). Now either eh ∈ E(G) or gh ∈ E(G) (else <B, g, e, C, h> = Kite). W.l.g. let eh ∈ E(G) ⇒ gh ∉ E(G) (else <e, g, h, B, u> = Kite). Similarly gf ∈ E(G) and ef ∉ E(G). Then D (C) is the only 4-vertex (3-vertex) of e (f) (else if j is another 4-vertex of e, then <B, u, f, e, j> = Chair). Similarly B (A) is the only 2-vertex (1-vertex) of g (h). Color C by 1, g by 2, B by 4, f by 3, u by 3, a contradiction.

Next define Q' = $\{X \in Q / XA \notin E(G)$ for some $A \in R\text{-}Q\}$. Then by **I**, $Q' \neq \phi$.

**Case 2.3: Cases 2.1 and 2.2 don't hold**.
Then $\forall$ Y $\in$ Q' and A $\in$ R-Q, YA $\notin$ E(G) (else let Z $\in$ Q' and B $\in$ R-Q be s.t. AZ, YB $\notin$ E(G). Then we get **Case 2.1 or 2.2)**.

Thus $\langle A, B, C, D\rangle = 2K_2$ whenever A, B $\in$ R-Q and C, D $\in$ Q'. **II**

W.l.g. let A, B $\in$ R-Q have color 1, 2 resply and C, D $\in$ Q' have color 3, 4 resply.

**Claim 1: 1-3 bi-color component N containing A and C is a path.**
Let P = {C, x, y, z,...., A} be a 3-1 bi-color path in N containing A and C. Now x is the unique 1-vertex of C (else if x' is another 1-vertex of C, then $\langle C, x, x', u, A\rangle$ = Chair). Also $deg_{\langle N\rangle}\, x = 2$ (else if y' is a third 2-vertex of x, then $\langle x, y, y', C, u\rangle$ = Chair). If N $\neq$ P, then $\exists$ w on P with $deg_{\langle N\rangle}\, w \geq 3$. Let w be the first such vertex on P. W.l.g. let w be colored 3. Let P $\supseteq$ {v, w', w, w''} and w''' be the third 1-vertex of w $\Rightarrow$ w $\notin$ {C, X} and $\langle w, w'', w''', w', v\rangle$ = chair, a contradiction. Hence N = P.

Now we proceed to show that $\exists$ a $K_r$ in G. Let N = {x $\in$ N(D) / x has same color as some A $\in$ R-Q}. Let a, b be a 1-vertex, 2-vertex of D resply. Then Ca, Cb $\notin$ E(G) (else $\langle a/b, C, D, u, A\rangle$ = Kite). Also a (b) is a unique 1-vertex (2-vertex) of D (else if g is another 1-vertex of D, then $\langle D, a, g, u, A\rangle$ = Chair). Then ab $\in$ E(G) (else if aB $\notin$ E(G), then $\langle D, a, b, u, B\rangle$ = Chair and if aB $\in$ E(G), then $\langle D, C, b, a, B\rangle$ = Chair). For every i-vertex X $\in$ R-Q, by **II**, DX $\notin$ E(G), and from above D has a unique i-vertex say x with xa, xb $\in$ E(G). Thus $\langle N\rangle$ is complete. **III**

Let y $\in$ Q-Q'. Also let x be an i-vertex in R-Q, and x' an i-vertex in N(D). Then yx' $\in$ E(G) (else $\langle x, u, y, D, x'\rangle$ = Kite). Thus $\langle (Q\text{-}Q')\cup N\rangle$ is complete. **IV**

**Claim 2: For every color i used in Q', D has another i-vertex.**
Let if possible Z be the only i-vertex of D where Z $\in$ Q'. Let T = {Z, g,...., A} be the i-1 path from Z to A. Now Dg $\notin$ E(G) (else $\langle g, D, Z, u, A\rangle$ = Kite). Also Za $\notin$ E(G) (else $\langle a, D, Z, u, A\rangle$ = Kite. Further a $\notin$ T (else $\langle a, a', a'', D, u\rangle$ = Chair, where T $\supseteq$ {a', a, a''}). Alter colors along T, color D by i and u by 4, a contradiction.

Let S = {x $\in$ N(D)-Q'/ x has same color as some X $\in$ Q'}.

**Claim 3: If s $\in$ S and m $\in$ N, then sm $\in$ E(G).**
Let V $\in$ R-Q be a j-vertex, W $\in$ Q' an i-vertex, m $\in$ N a j-vertex of D and s $\in$ S an i-vertex of D. Let if possible ms $\notin$ E(G). Then Wm $\notin$ E(G) (else $\langle m, D, W, u, V\rangle$ = Kite) and sV $\in$ E(G) (else $\langle D, m, s, u, V\rangle$ = Chair). But then $\langle D, m, W, s, V\rangle$ = Chair, a contradiction.

**Claim 4: $\langle S\rangle$ is complete.**
Let C, E $\in$ Q' and c, e $\in$ S be s.t. C, c; E, e have same colors.
Let if possible ce $\notin$ E(G). Let m $\in$ M have same color as some A $\in$ R-Q. By **Claim 3**, mc, me $\in$ E(G). Then Ac or Ae $\in$ E(G) (else $\langle D, c, e, u, A\rangle$ = Chair). W.l.g. let Ac $\in$ E(G). Then Ae $\in$ E(G) (else $\langle D, m, e, c, A\rangle$ = Kite). Further Ce $\in$ E(G) (else $\langle A, c, e, u, C\rangle$ = Chair). But then $\langle D, C, E, e, A\rangle$ = Kite, a contradiction.

Lastly let Z $\in$ Q-Q' and t $\in$ S. Also let T be an i-vertex in Q', and t an i-vertex in S. Let W $\in$ R-Q.

**Claim 5: Zt $\in$ E(G).**
Now Wt $\in$ E(G) (else $\langle Z, u, W, D, t\rangle$ = Kite). But then $\langle Z, u, T, W,t\rangle$ = Kite, a contradiction.

Thus from **Claim 3, 4, 5** and **II, III,** $\overline{N(D)}$ has a clique of size r > ω, a contradiction.

This proves **Theorem 2**.

**Corollary:** Reed's Conjecture is valid for every self-complementary graph without an induced Chair or Kite.

As Chair and Kite are complementary graphs, this follows from **Theorem 2**.

**Theorem 3:** If G is $\{K_2 \cup \overline{K_2}, H\}$-free then Reed's conjecture is valid for G.

Proof: Let G be a smallest $\{K_2 \cup \overline{K_2}, H\}$-free graph with $\chi(G) > \lceil \frac{\Delta+\omega+1}{2} \rceil$. Let $u \in V(G)$. By minimality, $\chi(G) - 1 \le \chi(G-u) \le \lceil \frac{\Delta(G-u)+\omega(G-u)+1}{2} \rceil \le \lceil \frac{\Delta+\omega+1}{2} \rceil < \chi(G)$. Thus $\chi(G-u) = \chi(G)-1 = \lceil \frac{\Delta+\omega+1}{2} \rceil$. Let C be a χ(G)-coloring of G, in which N(u) has r vertices with unique colors.

Then $\Delta \ge \deg u \ge r + 2(\lceil \frac{\Delta+\omega+1}{2} \rceil - r)$ and $r \ge \omega(G)+1$. **I**

Let $R = \{X \in N(u)/\ X$ receives a color unique in $N(u)\}$ and $S = \{A \in R/ \exists B \in R$ with $AB \notin E(G)\}$. As G is $K_2 \cup \overline{K_2}$-free, no two vertices in V(G)-N(u) have same color. Let $S' = \{a \in V(G)-N(u)/ \exists A \in S$ s.t. a, A have same color$\}$.

**Claim 1:** $<S'>$ is complete.

Let if possible $\exists\ a_1, a_2 \in S'$ with $A_1, A_2 \in S$ s.t. $a_1a_2 \notin E(G)$. Let $A_i$ have color i for i= 1, 2. Then $A_1a_2$, $A_2a_1 \in E(G)$ (else $<u, A_1, a_1, a_2> = K_2 \cup \overline{K_2}$) and $A_1A_2 \in E(G)$ (else $a_2$ ($A_1$) is the only 2-vertex (1-vertex) of $A_1$ ($a_2$). Color $A_1$ by 2, $a_2$ by 1, u by 1). Now $\exists\ A_3 \in S$ s.t. $A_1A_3 \notin E(G)$. Let $a_3 \in S'$. Clearly $A_1a_3$, $A_3a_1$, $a_1a_3 \in E(G) \Rightarrow a_2a_3 \in E(G)$ (else $<a_3, a_1, u, a_2> = K_2 \cup \overline{K_2}$). Again $a_3A_2 \in E(G)$ (else color $A_1$ by 3, $a_3$ by 2, $a_2$ by 1, u by 1). But then $<a_3, a_1, a_2, A_1, A_2, u> = H$, a contradiction.

Then $S \ne R$ (else $|S'| = r \ge \omega(G)+1$). Let $T = \{A \in R-S/ A$ has a unique color say α in G$\}$. Then $Aa_i \in E(G)\ \forall\ a_i \in S'$ (else let $A_j \in S$ s.t. $A_iA_j \notin E(G)$. Then $a_i$ is a unique i-vertex of $A_j$. Color $a_i$ by α, $A_j$ by i, u by j). Again as $<T \cup S'>$ is complete, ∃ color say 'δ' used in R but not in $T \cup S'$. Let B be a δ-vertex in R. Then as $B \notin T \cup S$, ∃ a δ-vertex $b \notin N(u)$. We show that $ba_i \in E(G)\ \forall\ a_i \in S'$. Let if possible $\exists\ a_i \in S'$ s.t. $ba_i \notin E(G)$. Let $A_i, A_j \in S$ and $a_i, a_j \in S'$ s.t. $A_iA_j \notin E(G)$ and $A_k$, $a_k$ have same colors for k = i, j. Than as before $A_ia_j$, $A_ja_i$, $a_ia_j \in E(G)$. Also $A_ib$, $Ba_i$, $AB \in E(G)$. As $B \notin S$, $BA_i$, $BA_j \in E(G)$. Also $ba_j \in E(G)$ (else $<a_i, a_j, b, u> = K_2 \cup \overline{K_2}$) and $Ba_j \in E(G)$ (else b is the only δ-vertex of $a_j$. Color $a_j$ by δ, b by i, $A_i$ by j, u by i). But then $<a_j, a_i, b, A_i, B, u> = H$, a contradiction. Hence $ba_i \in E(G)\ \forall\ a_i \in S'$. Let $W = \{b \notin N(u)/ \exists\ B \notin R-(T \cup S)$ s.t. b, B have same color$\}$. Then $<W \cup T \cup S'> = K_r$, a contradiction as $r > \omega(G)$.

This proves **Theorem 3.**

**Theorem 4:** If G is $\{2K_2, M\}$-free, then Reed's conjecture is valid for G.

Proof: Let G be a smallest $\{K_2 \cup \overline{K_2}, M\}$-free graph with $\chi(G) > \lceil \frac{\Delta+\omega+1}{2} \rceil$. Let $u \in V(G)$. By minimality, $\chi(G) - 1 \le \chi(G-u) \le \lceil \frac{\Delta(G-u)+\omega(G-u)+1}{2} \rceil \le \lceil \frac{\Delta+\omega+1}{2} \rceil < \chi(G)$. Thus $\chi(G-u) =$ $\chi(G)-1 = \lceil \frac{\Delta+\omega+1}{2} \rceil$. Let C be a $\chi(G)$-coloring of G, in which N(u) has r vertices with unique colors in N(u). Then $\Delta \geq \deg u \geq r + 2(\lceil \frac{\Delta+\omega+1}{2} \rceil - r)$ and $r \geq \omega(G)+1$. **I**

Let $R = \{X \in N(u) /\ X$ receives a unique color in $N(u)\}$ and $S = \{A_i \in R /\ \exists\ A_j \in R$ with $A_iA_j \notin E(G)\}$. Then $|R| = r$.

**Claim 1:** $\forall\ A_i \in S, \exists\ a_i \in N(u)-R$ s.t. $a_iA_j \in E(G)$ for $1 \leq j \neq i \leq r$.
Let if possible $\exists\ A_i \in S$ s.t. no $a_i \in N(u)-R$ is adjacent to all $A_k$, $k \neq i$. Now let $A_j \in S$ s.t. $A_iA_j \notin E(G)$. Then $A_j$ has an i-vertex say $x \in N(u)-R$. By assumption $xA_m \notin E(G)$ for some m. As G is $2K_2$-free, x has no k–vertex in G. Color x by k. Thus all i-vertices of $A_j$ can be colored by colors other than i and $A_j$ has no i-vertex. Color $A_j$ by i and u by j, a contradiction. Hence **Claim 1** holds.

Let $S' = \{a_i \in N(u)-R / A_i \in S$ and $a_iA_j \in E(G)$ for $1 \leq j \neq i \leq r, A_j \in S\}$.

**Claim 2:** $<S'>$ is complete.
Let if possible $\exists\ a_1, a_2 \in S'$ with $A_1, A_2 \in S$ s.t. $a_1a_2 \notin E(G)$. Then $A_1A_2 \in E(G)$ (else by **Claim 1**, $<A_1, A_2, a_1, a_2> = 2K_2$) and $\exists\ A_3 \in S$ s.t. $A_1A_3 \notin E(G)$. Now let $a_3 \in S'$. Then by **Claim 1**, $a_1A_3$, $A_1a_3 \in E(G)$ and $a_1a_3 \in E(G)$ (else $<A_1, a_3, A_3, a>_1 = 2K_2$) $\Rightarrow a_3a_2 \notin E(G)$ (else $<u, A_1, a_3, a_1, A_3, a_2> = M$) $\Rightarrow A_3A_2 \in E(G)$ (else by **Claim 1**, $<A_3, a_2, A_2, a_3> = 2K_2$). As $A_2 \in S$, $\exists\ A_4 \in S$ s.t. $A_2A_4 \notin E(G)$. Let $a_4 \in S'$. Then $a_4a_2 \in E(G)$. As G is $2K_2$-free either $a_4a_1 \in E(G)$ or $a_4a_3 \in E(G)$. But then $<u, A_2, a_4, a_2, A_4, a_1/a_3> = M$, a contradiction.
Hence **Claim 2** holds.

If $R = S$, then by **I**, $\omega(G) \geq |S'| = r \geq \omega(G)+1$, a contradiction. If $R \supset S$, then for every $A_i \in R-S$, $A_iA_j \in E(G)$ where $A_j \in S$ and $A_ia_k \in E(G)\ \forall\ a_k \in S'$ by **Claim 1**. Thus $<(R-S)\cup S'>$ is complete and $\omega(G) \geq r \geq \omega(G)+1$, a contradiction.

This proves **Theorem 4**.